\documentclass[reqno]{article}%
\usepackage[intlimits]{amsmath}
\usepackage{graphicx}
\usepackage{amsfonts}
\usepackage{amssymb}%
\newtheorem{theorem}{Theorem}

\newtheorem{lemma}[theorem]{Lemma}

\newtheorem{problem}[theorem]{Problem}

\newenvironment{proof}[1][Proof]{\noindent{\textbf {#1}  }}  {\hfill$\Box$}

\begin{document}

\title{Edge distribution of graphs with few induced copies of a given graph}
\author{V. Nikiforov\ \\{\small Department of Mathematical Sciences}\\{\small University of Memphis, Memphis, TN 38152}\\{\small e-mail:} {\small vnikifrv@memphis.edu}}
\date{November 26, 2003}
\maketitle

\begin{abstract}
We show that if a simple graph contains few induced copies of a given graph,
then its edges are distributed rather unevenly.\ 

More precisely, for all $\varepsilon>0$ and $r\geq2,$ there exist $\xi
=\xi\left(  \varepsilon,r\right)  >0$ and $L=L\left(  \varepsilon,r\right)  $
such that, for every graph $H$ of order $r,$ and every graph $G$ of
sufficiently large order $n,$ the following assertion holds.

If $G$ contains fewer than $\xi n^{r}$ copies of $H$, then there exists a
partition $V\left(  G\right)  =\cup_{i=0}^{q}V_{i}$ with $\left\vert
V_{0}\right\vert <q\leq L,$ such that $\left\vert V_{i}\right\vert
=\left\lfloor n/q\right\rfloor ,$ and%
\[
e\left(  V_{i}\right)  <\varepsilon\binom{\left\vert V_{i}\right\vert }%
{2}\ \ \ \text{or}\ \ \ e\left(  V_{i}\right)  >\left(  1-\varepsilon\right)
\binom{\left\vert V_{i}\right\vert }{2}%
\]
for every $i\in\left[  q\right]  .$

In particular, for all $\varepsilon>0$ and $r\geq2,$ there exist $\xi
=\xi\left(  \varepsilon,r\right)  >0$ and $L=L\left(  \varepsilon,r\right)  $
such that, for every graph $G$ of sufficiently large order $n,$ the following
assertion holds.

If $G$ has fewer than $\xi n^{r}$ $r$-cliques, then there exists a partition
$V\left(  G\right)  =\cup_{i=0}^{q}V_{i}$ with $\left\vert V_{0}\right\vert
<q\leq L$ such that
\[
\left\vert V_{i}\right\vert =\left\lfloor n/q\right\rfloor ,\text{ \ \ and
\ \ }e\left(  W_{i}\right)  <\varepsilon\binom{\left\vert V_{i}\right\vert
}{2}%
\]
for every $i\in\left[  q\right]  .$

We derive also a number of related results.

\end{abstract}

\section{Introduction}

Our graph-theoretic notation is standard (e.g., see \cite{Bol98}); thus we
write $G\left(  n,m\right)  $ for a graph of order $n$ and size $m$. Given two
graphs $H$ and $G$ we write $k_{H}\left(  G\right)  $ for the number of
induced copies of $H$ in $G;$ $k_{r}\left(  G\right)  $ stands for $k_{K_{r}%
}\left(  G\right)  .$ If $U\subset V\left(  G\right)  ,$ we write $e\left(
U\right)  $ for $e\left(  G\left[  U\right]  \right)  ,$ $k_{H}\left(
U\right)  $ for $k_{H}\left(  G\left[  U\right]  \right)  ,$ and $k_{r}\left(
U\right)  $ for $k_{r}\left(  G\left[  U\right]  \right)  $. A partition
$V=\cup_{i=0}^{k}V_{i}$ is called \emph{equitable,} if $\left\vert
V_{0}\right\vert <k,$ and $\left\vert V_{1}\right\vert =...=\left\vert
V_{k}\right\vert .$ A set of cardinality $k$ is called a $k$-set.

In \cite{Erd79} Erd\H{o}s raised the following problem (see also \cite{Bol78},
p. 363).

\begin{problem}
\label{Erp} Let $c>0.$ Suppose $G=G(n,\left\lfloor cn^{2}\right\rfloor )$ is
such that
\[
e\left(  W\right)  \geq\left(  c/4+o\left(  1\right)  \right)  n^{2}%
\]
for every $W\subset V\left(  G\right)  $ with $\left\vert W\right\vert
=\left\lfloor n/2\right\rfloor .$ Then, for every fixed $r$ and sufficiently
large $n,$ the graph $G$ contains $K_{r}$.
\end{problem}

This problem was solved recently in \cite{Nik01}, where the following more
general result was proved.

\begin{theorem}
\label{Ers} For every $c>0$ and $r\geq3,$ there exists $\beta=\beta(c,r)>0$
such that, for every $K_{r}$-free graph $G=G(n,m)$ with $m\geq cn^{2}$, there
exists a partition $V\left(  G\right)  =V_{1}\cup V_{2}$ with $\left\vert
V_{1}\right\vert =\left\lfloor n/2\right\rfloor ,$ $\left\vert V_{2}%
\right\vert =\left\lceil n/2\right\rceil ,$ and
\begin{equation}
e\left(  V_{1},V_{2}\right)  >\left(  1/2+\beta\right)  m. \label{maxcutlb}%
\end{equation}

\end{theorem}

In fact, (\ref{maxcutlb}) is a lower bound on the MaxCut function for dense
$K_{r}$-free graphs; note, that it differs significantly from those found in
\cite{Alo96}, \cite{ABKS03} and \cite{AKS03}. We obtain a similar result about
judicious partitions in Theorem \ref{Ers2}.

Kohayakawa and R\"{o}dl \cite{KoRo03} gave another solution to Problem
\ref{Erp}; however, their method does not imply Theorem \ref{Ers}.

One of our goals in this note is to extend Theorem \ref{Ers}. We first prove
the following basic result.

\begin{theorem}
\label{maint} For all $\varepsilon>0$ and $r\geq2,$ there exist $\xi
=\xi\left(  \varepsilon,r\right)  >0$ and $L=L\left(  \varepsilon,r\right)  $
such that, for every graph $G$ of sufficiently large order $n,$ the following
assertion holds.

If $k_{r}\left(  G\right)  <\xi n^{r}$, then there exists an equitable
partition $V\left(  G\right)  =\cup_{i=0}^{q}V_{i}$ with $q<L,$ and
\[
e\left(  V_{i}\right)  <\varepsilon\binom{\left\vert V_{i}\right\vert }{2}%
\]
for every $i\in\left[  q\right]  .$
\end{theorem}

From this assertion we shall deduce that the conclusion of Theorem \ref{Ers}
remains essentially true under considerably weaker stipulations.

\begin{theorem}
\label{Ers1} For all $c>0$ and $r\geq3$, there exist $\xi=\xi\left(
c,r\right)  >0$ and $\beta=\beta(c,r)>0$such that, for $n$ sufficiently large
and every graph $G=G(n,m)$ with $m\geq cn^{2}$, the following assertion holds.

If $k_{r}\left(  G\right)  <\xi n^{r}$, then there exists a partition
$V\left(  G\right)  =V_{1}\cup V_{2}$ with $\left\vert V_{1}\right\vert
=\left\lfloor n/2\right\rfloor ,$ $\left\vert V_{2}\right\vert =\left\lceil
n/2\right\rceil ,$ and
\[
e\left(  V_{1},V_{2}\right)  >\left(  1/2+\beta\right)  m.
\]

\end{theorem}

We deduce also a number of related results, in particular, the following
analogue of Theorem \ref{maint}.

\begin{theorem}
\label{maintx} For all $\varepsilon>0$ and $r\geq2,$ there exist $\xi
=\xi\left(  \varepsilon,r\right)  >0$ and $L=L\left(  \varepsilon,r\right)  $
such that, for every graph $H$ of order $r,$ and every graph $G$ of
sufficiently large order $n,$ the following assertion holds.

If $k_{H}\left(  G\right)  <\xi n^{r}$, then there exists an equitable
partition $V\left(  G\right)  =\cup_{i=0}^{q}V_{i}$ with $q<L$ such that
\[
e\left(  V_{i}\right)  <\varepsilon\binom{\left\vert V_{i}\right\vert }%
{2}\ \ \ \text{or}\ \ \ e\left(  V_{i}\right)  >\left(  1-\varepsilon\right)
\binom{\left\vert V_{i}\right\vert }{2}%
\]
for every $i\in\left[  q\right]  .$
\end{theorem}

Observe that, although Theorem \ref{maintx} is a fairly general result, it
does not imply Theorem \ref{maint} or its counterpart for independent $r$-sets.

Finally, we prove the following assertion that looks likely to be useful in
Ramsey type applications; we shall investigate this topic in a forthcoming note.

\begin{theorem}
\label{Rams} For all $\varepsilon>0,$ $r\geq2$ and $k\geq2,$ there exist
$\delta=\delta\left(  \varepsilon,r\right)  >0,$ $\xi=\xi\left(
\varepsilon,r\right)  >0$\ and $L=L\left(  \varepsilon,r,k\right)  $ such
that,\ for every graph $G$ of sufficiently large order $n,$ the following
assertion holds.

If $V\left(  G\right)  =\cup_{i=0}^{k}V_{i}$ is a $\delta$-uniform partition
such that%
\[
k_{r}\left(  V_{i}\right)  \leq\xi\left\vert V_{i}\right\vert ^{r}\text{
\ \ or \ \ }k_{r}\left(  \overline{G\left[  V_{i}\right]  }\right)  \leq
\xi\left\vert V_{i}\right\vert ^{r}%
\]
for every $i\in\left[  k\right]  ,$ then there exists an $\varepsilon$-uniform
partition $V\left(  G\right)  =\cup_{i=0}^{q}W_{i}$ with $k\leq q\leq L$ such
that
\[
e\left(  W_{i}\right)  <\varepsilon\binom{\left\vert W_{i}\right\vert }%
{2}\text{ \ \ or \ \ }e\left(  W_{i}\right)  >\left(  1-\varepsilon\right)
\binom{\left\vert W_{i}\right\vert }{2}%
\]
for every $i\in\left[  q\right]  .$
\end{theorem}

The rest of the note is organized as follows. First we introduce some
additional notation, then we prove Theorem \ref{maint} in Section
\ref{Mproof}, extend it in Section \ref{Cor}, and use it in Section
\ref{Bipar} to prove Theorem \ref{Ers1}. In Section \ref{Spar} we prove
Theorem \ref{Rams} and, finally, in Section \ref{Mxproof} we prove Theorem
\ref{maintx}.

A few words about our proofs seem necessary. We apply continually
Szemer\'{e}di's uniformity lemma (SUL) in a rather routine manner. However,
for the reader's sake, we always provide the necessary details, despite repetitions.

\subsection{Notation}

Suppose $G$ is a graph. For a vertex $u\in V\left(  G\right)  ,$ we write
$\Gamma\left(  u\right)  $ for the set of vertices adjacent to $u$. If
$A,B\subset V\left(  G\right)  $ are nonempty disjoint sets, we write
$e\left(  A,B\right)  $ for the number of $A-B$ edges and set%
\[
d\left(  A,B\right)  =\frac{e\left(  A,B\right)  }{\left\vert A\right\vert
\left\vert B\right\vert }.
\]

Given a partition $V=\cup_{i=0}^{k}V_{i}$,\ we occasionally call the sets
$V_{1},...,V_{k}$ \emph{clusters} of the partition.

For general notions and definitions related to Szemer\'{e}di's uniformity
lemma (SUL), see, e.g. \cite{KoSi93}, or \cite{Bol98}. In our exposition we
shall systematically replace \textquotedblleft regularity\textquotedblright%
\ by \textquotedblleft uniformity\textquotedblright, thus \textquotedblleft%
$\varepsilon$-uniform\textquotedblright\ will stand for \textquotedblleft%
$\varepsilon$-regular\textquotedblright.

Let $\varepsilon>0$. A partition $V\left(  G\right)  =\cup_{i=0}^{k}V_{i}$ is
called $\varepsilon$\emph{-uniform,} if it is equitable, and at most
$\varepsilon k^{2}$ pairs $\left(  V_{i},V_{j}\right)  $ are not $\varepsilon$-uniform.

\section{\label{Mproof}Proof of Theorem \ref{maint}}

In our proof of Theorem \ref{maint} and later we shall use SUL in the
following form.

\begin{theorem}
[Szemer\'{e}di's Uniformity Lemma]\label{SUL} Let $l\geq1$, $\varepsilon>0$.
There exists $M=M\left(  \varepsilon,l\right)  $ such that, for every graph
$G$ of sufficiently large order, there exists an $\varepsilon$-uniform
partition $V\left(  G\right)  =\cup_{i=0}^{k}V_{i}$ with $l\leq k\leq M.$
\end{theorem}

In addition, we need the following basic properties of $\varepsilon$-uniform
pairs (see \cite{KoSi93}, Facts 1.4 and 1.5.)

\begin{lemma}
\label{Xsec}Let $\varepsilon>0,$ $r\geq1,$ and $\left(  A,B\right)  $ be an
$\varepsilon$-uniform pair with $d\left(  A,B\right)  =d.$ If $Y\subset B$ and
$\left(  d-\varepsilon\right)  ^{r-1}\left\vert Y\right\vert >\varepsilon
\left\vert B\right\vert ,$ then there are at most $\varepsilon r\left\vert
A\right\vert ^{r}$ $r$-sets $R\subset A$ such that%
\[
\left\vert \left(  \cap_{u\in R}\Gamma\left(  u\right)  \right)  \cap
Y\right\vert \leq\left(  d-\varepsilon\right)  ^{r}\left\vert Y\right\vert .
\]

\end{lemma}

\begin{lemma}
\label{Slicel} Let $0<\varepsilon<\alpha,$ and let $\left(  A,B\right)  $ be
an $\varepsilon$-uniform pair. If $A^{\prime}\subset A,$ $B^{\prime}\subset B$
and $\left\vert A^{\prime}\right\vert \geq\alpha\left\vert A\right\vert ,$
$\left\vert B^{\prime}\right\vert \geq\alpha\left\vert B\right\vert ,$ then
$\left(  A^{\prime},B^{\prime}\right)  $ is an $\varepsilon^{\prime}$-uniform
pair with $\varepsilon^{\prime}=\max\left\{  \varepsilon/\alpha,2\varepsilon
\right\}  .$
\end{lemma}

It is straightforward to deduce the following assertion from Lemma \ref{Xsec}.

\begin{lemma}
\label{XPle} Let $r\geq1,$ $0<2\varepsilon^{1/r}<d\leq1,$ and let $\left(
A,B\right)  $ be an $\varepsilon$-uniform pair with $d\left(  A,B\right)  =d.$
There are at most $\varepsilon r\left\vert A\right\vert ^{r}$ $r$-sets
$R\subset A$ such that%
\[
\left\vert \left(  \cap_{u\in R}\Gamma\left(  u\right)  \right)  \cap
B\right\vert \leq\varepsilon\left\vert B\right\vert .
\]

\end{lemma}

The following simple lemma will play a crucial role in our proofs.

\begin{lemma}
[Scooping Lemma]\label{partL} Let $\varepsilon>0,$ and let $s$ be integer with
$0<s\leq\varepsilon n.$ For every graph $G$ of order $n,$ if $e\left(
G\right)  \leq\varepsilon^{3}\binom{n}{2},$ then there exists a partition
$V\left(  G\right)  =\cup_{i=0}^{k}V_{i}$ such that $\left\vert V_{0}%
\right\vert \leq\left\lceil \varepsilon n\right\rceil ,$ and
\[
\left\vert V_{i}\right\vert =s,\text{ \ \ \ }e\left(  V_{i}\right)
<\varepsilon\binom{s}{2}%
\]
for every $i\in\left[  k\right]  .$
\end{lemma}

\begin{proof}
Select a sequence of sets $V_{1},...,V_{k}$ as follows: select $V_{1}$ by
\[
e\left(  V_{1}\right)  =\min\left\{  e\left(  U\right)  :U\subset V\left(
G\right)  ,\text{ }\left\vert U\right\vert =s\right\}  ;
\]
having selected $V_{1},...,V_{i},$ if $n-is\leq\left\lceil \varepsilon
n\right\rceil $ stop the sequence, else select $V_{i+1}$ by%
\[
e\left(  V_{i+1}\right)  =\min\left\{  e\left(  U\right)  :U\subset V\left(
G\right)  \backslash\left(  \cup_{j=1}^{i}V_{j}\right)  ,\text{ }\left\vert
U\right\vert =s\right\}  .
\]

Let $V_{k}$ be the last selected set; set $V_{0}=V\left(  G\right)
\backslash\left(  \cup_{i=1}^{k}V_{i}\right)  .$ The stop condition implies
$\left\vert V_{0}\right\vert \leq\left\lceil \varepsilon n\right\rceil .$ For
every $i\in\left[  k\right]  ,$ the way we choose $V_{i}$ implies%
\begin{align*}
e\left(  V_{i}\right)   &  \leq\frac{e\left(  G\right)  }{\binom{n-\left(
i-1\right)  s}{2}}\binom{s}{2}\leq\frac{\varepsilon^{3}n\left(  n-1\right)
}{\left(  n-\left(  i-1\right)  s\right)  \left(  n-\left(  i-1\right)
s-1\right)  }\binom{s}{2}\\
&  <\frac{\varepsilon^{3}n\left(  n-1\right)  }{\left(  \left\lceil
\varepsilon n\right\rceil +1\right)  \left\lceil \varepsilon n\right\rceil
}\binom{s}{2}<\varepsilon\binom{s}{2},
\end{align*}
so the partition $V\left(  G\right)  =\cup_{i=0}^{k}V_{i}$ has the required properties.
\end{proof}

\begin{proof}
[Proof of Theorem \ref{maint}]Setting $q=L\left(  \varepsilon,2\right)  =1,$
$\xi\left(  \varepsilon,2\right)  =\varepsilon,$ the theorem holds trivially
for $r=2.$ To prove it for $r>2$ we apply induction on $r$ - assuming it holds
for $r,$ we shall prove it for $r+1.$

Observe that it suffices to find $\xi=\xi\left(  \varepsilon,r+1\right)  >0$
and $L=L\left(  \varepsilon,r+1\right)  $ such that, if $G$ is a graph of
sufficiently large order $n,$ and $k_{r+1}\left(  G\right)  <\xi n^{r+1},$
then there exists a partition $V\left(  G\right)  =\cup_{i=0}^{q}W_{i}$ such that:

\emph{(i)} $q\leq L;$\ 

\emph{(ii)} $\left\vert W_{0}\right\vert <6\varepsilon n,$ $\left\vert
W_{1}\right\vert =...=\left\vert W_{q}\right\vert ;$

\emph{(iii)} for every $i\in\left[  q\right]  ,$ $e\left(  W_{i}\right)
<\varepsilon\binom{\left\vert W_{i}\right\vert }{2}.$

Indeed, distributing evenly among the sets $W_{1},...,W_{q}$ as many as
possible of the vertices of $W_{0},$ we obtain a partition $V\left(  G\right)
=\cup_{i=0}^{q}V_{i}$ with $\left\vert V_{0}\right\vert <q,$ and
\[
\left\vert V_{i}\right\vert =\left\lfloor n/q\right\rfloor ,\ \text{\ \ \ \ }%
e\left(  V_{i}\right)  <2\varepsilon\binom{\left\vert V_{i}\right\vert }{2},
\]
for every $i\in\left[  q\right]  ,$ as required.

For convenience we shall outline first our proof. For $\delta$ appropriately
small, applying SUL, we find a $\delta$-uniform partition $V\left(  G\right)
=\cup_{i=0}^{k}V_{i}.$ Note that, if $k_{r}\left(  V_{i}\right)  $ is
proportional to $n^{r}$, and $V_{i}$ is incident to a substantially dense
$\delta$-uniform pair, then, by Lemma \ref{XPle}, there are substantially many
$\left(  r+1\right)  $-cliques in $G.$ Therefore, for every $V_{i},$ either
$k_{r}\left(  V_{i}\right)  $ is small or $V_{i}$ is essentially isolated.

Let $V^{\prime\prime}$ be the union of the essentially isolated clusters; set
$V^{\prime}=V\backslash\left(  V^{\prime\prime}\cup V_{0}\right)  $.

\emph{1 Partitioning of }$V^{\prime}$

By the induction hypothesis, we partition each nonisolated $V_{i}$ into a
bounded number of sparse sets $Y_{ij}$ and a small exceptional set; the
exceptional sets are collected in $X^{\prime}.$ Note that, although the sets
$Y_{ij}$ are sparse, they are not good for our purposes, for their cardinality
may vary with $i$. To overcome this obstacle, we first select a sufficiently
small integer $s$ proportional to $n.$ Then, by the Scooping Lemma, we
partition each of the sets $Y_{ij}$ into sparse sets of cardinality exactly
$s$ and a small exceptional set; the exceptional sets are added to $X^{\prime
}$.

\emph{2 Partitioning of }$V^{\prime\prime}$

We partition $V^{\prime\prime}$ into sparse sets of size $s$ and a small
exceptional set $X_{0}$. If $\left\vert V^{\prime\prime}\right\vert $ is
small, we set $X_{0}=V^{\prime\prime},$ and complete the partition. If
$\left\vert V^{\prime\prime}\right\vert $ is substantial, then $G\left[
V^{\prime\prime}\right]  $ must be sparse, for it consists of essentially
isolated clusters. Applying the Scooping Lemma to $G\left[  V^{\prime\prime
}\right]  $, we partition $V^{\prime\prime}$ into sparse sets of cardinality
$s$ and a small exceptional set $X_{0}$.

Let $W_{1},...,W_{q}$ be the sets of cardinality $s$ obtained during the
partitioning of $V^{\prime}$ and $V^{\prime\prime}$. Set $W_{0}=V_{0}\cup
X_{0}\cup X^{\prime};$ the choice of $\delta$ implies $\left\vert
W_{0}\right\vert <6\varepsilon n,$ so the partition $V\left(  G\right)
=\cup_{i=0}^{q}W_{i}$ satisfies \emph{(i)-(iii)}.

Let us now give the details. Assume $\varepsilon$ sufficiently small and set%
\begin{align}
l  &  =\max\left\{  \left\lceil \frac{1}{\varepsilon^{5}}\right\rceil
,\frac{1}{4\varepsilon L\left(  \varepsilon^{3},r\right)  }\right\}
,\label{defl}\\
\delta &  =\min\left\{  \frac{\xi\left(  \varepsilon^{3},r\right)  }%
{r+1},\frac{\varepsilon^{5r}}{16^{r}}\right\}  .\label{defdel}\\
L  &  =L\left(  \varepsilon,r+1\right)  =\frac{8M\left(  \delta,l\right)
L\left(  \varepsilon^{3},r\right)  }{\varepsilon}.\label{defL}\\
\xi &  =\xi\left(  \varepsilon,r+1\right)  =\frac{\delta^{2}}{\left(
2M\left(  \delta,l\right)  \right)  ^{r+1}} \label{defksi}%
\end{align}

Let $G$ be a graph of sufficiently large order $n,$ and let $k_{r+1}\left(
G\right)  <\xi n^{r+1}.$ Applying SUL, we find a $\delta$-uniform partition
$V\left(  G\right)  =\cup_{i=0}^{k}V_{i}$ with $l\leq k\leq M\left(
\delta,l\right)  $. Set $t=\left\vert V_{1}\right\vert $ and observe that
\begin{equation}
\frac{n}{2k}\leq\left(  1-\delta\right)  \frac{n}{k}<t\leq\frac{n}{k}.
\label{bndt}%
\end{equation}
Assume that there exist a cluster $V_{i}$ with $k_{r}\left(  V_{i}\right)
>\xi\left(  \varepsilon^{3},r\right)  t^{r},$ and a $\delta$-uniform pair
$\left(  V_{i},V_{j}\right)  $ with
\[
d\left(  V_{i},V_{j}\right)  >2\delta^{1/r}.
\]
Applying Lemma \ref{XPle} with $A=V_{i}$ and $B=V_{j},$ we find that there are
at least
\[
\xi\left(  \varepsilon^{3},r\right)  t^{r}-\delta rt^{r}\geq\delta t^{r}%
\]
$r$-cliques $R\subset V_{i}$ such that%
\[
\left\vert \left(  \cap_{u\in R}\Gamma\left(  u\right)  \right)  \cap
V_{j}\right\vert >\delta t.
\]
Hence, there are at least $\delta^{2}t^{r+1}$ $\left(  r+1\right)  $-cliques
inducing an $r$-clique in $V_{i}$ and a vertex in $V_{j}$. Therefore, from
(\ref{bndt}) and (\ref{defksi}), we find that
\begin{align*}
k_{r+1}\left(  G\right)   &  \geq\delta^{2}t^{r+1}>\delta^{2}\left(
\frac{1-\delta}{k}\right)  ^{r+1}n^{r+1}\\
&  >\frac{\delta^{2}}{\left(  2M\left(  \delta,l\right)  \right)  ^{r+1}%
}n^{r+1}=\xi\left(  \varepsilon,r+1\right)  n^{r+1},
\end{align*}
a contradiction. Therefore, if $k_{r}\left(  V_{i}\right)  \geq\xi\left(
\varepsilon^{3},r\right)  t^{r},$ then every $\varepsilon$-uniform pair
$\left(  V_{i},V_{j}\right)  $ satisfies $d\left(  V_{i},V_{j}\right)
\leq2\delta^{1/r}.$ Let
\[
I^{\prime}=\left\{  i:i\in\left[  k\right]  \text{, \ }k_{r}\left(
V_{i}\right)  \leq\xi\left(  \varepsilon^{3},r\right)  t^{r}\right\}  ,\text{
\ \ }I^{\prime\prime}=\left[  k\right]  \backslash I^{\prime}.
\]

First we shall partition $V^{\prime}=\cup_{i\in I^{\prime}}V_{i}.$ Set%
\begin{equation}
s=\left\lfloor \frac{\varepsilon n}{4kL\left(  \varepsilon^{3},r\right)
}\right\rfloor \label{defs}%
\end{equation}
and observe that
\begin{equation}
\frac{n}{s}<\frac{8kL\left(  \varepsilon^{3},r\right)  }{\varepsilon}\leq
\frac{8M\left(  \delta,l\right)  L\left(  \varepsilon^{3},r\right)
}{\varepsilon}=L\left(  \varepsilon,r+1\right)  . \label{props}%
\end{equation}
For every $i\in I^{\prime},$ by the induction hypothesis, we find an equitable
partition $V_{i}=\cup_{j=0}^{m_{i}}Y_{ij}$ with $\left\vert Y_{i0}\right\vert
<m_{i}\leq L\left(  \varepsilon^{3},r\right)  ,$ and\
\[
e\left(  Y_{ij}\right)  \leq\varepsilon^{3}\binom{\left\vert Y_{ij}\right\vert
}{2}%
\]
for every $j\in\left[  m_{i}\right]  .$ Also, for every $i\in I^{\prime}$ and
$j\in\left[  m_{i}\right]  ,$ (\ref{defs}) and (\ref{bndt}) imply%
\[
s\leq\frac{\varepsilon n}{4kL\left(  \varepsilon^{3},r\right)  }%
<\varepsilon\frac{t}{2L\left(  \varepsilon^{3},r\right)  }\leq\varepsilon
\frac{t}{2m_{i}}\leq\varepsilon\left\lfloor \frac{t}{m_{i}}\right\rfloor
=\varepsilon\left\vert Y_{ij}\right\vert .
\]
Hence, we apply the Scooping Lemma to the graph $G\left[  Y_{ij}\right]  ,$
and find a partition $Y_{ij}=\cup_{q=0}^{p_{ij}}W_{ijq}$ with $\left\vert
W_{ij0}\right\vert \leq\left\lceil \varepsilon\left\vert Y_{ij}\right\vert
\right\rceil $ such that
\[
\left\vert W_{ijq}\right\vert =s\text{ \ \ and \ \ }e\left(  W_{ijq}\right)
<\varepsilon\binom{s}{2}%
\]
for every $q\in\left[  p_{ij}\right]  .$ Setting $X^{\prime}=\left(
\cup_{i\in I}Y_{i0}\right)  \cup\left(  \cup_{i\in I}\cup_{j=1}^{m_{i}}%
W_{ij0}\right)  ,$ we obtain
\begin{align}
\left\vert X^{\prime}\right\vert  &  =\left\vert \left(  \cup_{i\in I}%
Y_{i0}\right)  \right\vert +\left\vert \left(  \cup_{i\in I}\cup_{j=1}^{m_{i}%
}W_{ij0}\right)  \right\vert <\sum_{i\in I^{\prime}}m_{i}+2\varepsilon
\sum_{i\in I^{\prime}}\sum_{j=1}^{m_{i}}\left\vert Y_{ij}\right\vert
\nonumber\\
&  <kL\left(  \varepsilon^{3},r\right)  +2\varepsilon\sum_{i\in I^{\prime}%
}m_{i}\left\lfloor \frac{t}{m_{i}}\right\rfloor \leq kL\left(  \varepsilon
^{3},r\right)  +2\varepsilon n<3\varepsilon n. \label{x1up}%
\end{align}

Denote by $h$ the number of the sets $W_{ijq}$ $\left(  i\in I^{\prime},\text{
}j\in\left[  m_{i}\right]  ,\text{ }q\in\left[  p_{ij}\right]  \right)  ,$ and
renumber them sequentially from $1$ to $h$.\ So far we have a partition
$V^{\prime}=X^{\prime}\cup\left(  \cup_{i=1}^{h}W_{i}\right)  $ with
$\left\vert X^{\prime}\right\vert <3\varepsilon n$ such that%
\[
\left\vert W_{i}\right\vert =s,\text{ \ \ \ }e\left(  W_{i}\right)
<\varepsilon\binom{s}{2}%
\]
for every $i\in\left[  h\right]  .$

Next we shall partition the set $V^{\prime\prime}=\cup_{i\in I^{\prime\prime}%
}V_{i}.$ We may assume that $\left\vert V^{\prime\prime}\right\vert
\geq\varepsilon n,$ else, setting $W_{0}=V_{0}\cup X^{\prime}\cup
V^{\prime\prime}$, from (\ref{x1up}), we have $W_{0}<5\varepsilon n,$\ and, in
view of (\ref{props}), the proof is completed. Obviously,%
\begin{equation}
e\left(  V^{\prime\prime}\right)  =\sum_{i\in I^{\prime\prime}}e\left(
V_{i}\right)  +\sum_{i,j\in I^{\prime\prime},\text{ }i<j}e\left(  V_{i}%
,V_{j}\right)  \leq k\binom{t}{2}+e_{1}+e_{2}, \label{ev}%
\end{equation}
where,
\begin{align*}
e_{1}  &  =\sum_{i,j\in I^{\prime\prime},\text{ }i<j}\left\{  e\left(
V_{i},V_{j}\right)  :\left(  V_{i},V_{i}\right)  \text{ is }\delta
\text{-uniform}\right\}  ,\\
e_{2}  &  =\sum_{i,j\in I^{\prime\prime},\text{ }i<j}\left\{  e\left(
V_{i},V_{j}\right)  :\left(  V_{i},V_{i}\right)  \text{ is not }%
\delta\text{-uniform}\right\}
\end{align*}
Since at most $\delta k^{2}$ pairs $\left(  V_{i},V_{i}\right)  ,$ $\left(
1\leq i<j\leq k\right)  ,$ are not $\delta$-uniform, it follows that
\begin{equation}
e_{2}\leq\delta k^{2}t^{2}\leq\delta n^{2}. \label{s2up}%
\end{equation}
Recall that if $i,j\in I^{\prime\prime},$ $i<j,$ and the pair $\left(
V_{i},V_{j}\right)  $ is $\delta$-uniform, then $d\left(  V_{i},V_{j}\right)
<2\delta^{1/r}.$ Therefore,%
\[
e_{1}\leq\binom{k}{2}2\delta^{1/r}t^{2}\leq\delta^{1/r}k^{2}t^{2}\leq
\delta^{1/r}n^{2}.
\]
Hence, (\ref{ev}), (\ref{s2up}), (\ref{bndt}), and (\ref{defl}) imply
\begin{align*}
e\left(  V^{\prime\prime}\right)   &  \leq\frac{kt^{2}}{2}+\delta n^{2}%
+\delta^{1/r}n^{2}\leq\left(  \frac{1}{2k}+2\delta^{1/r}\right)  n^{2}%
\leq\left(  \frac{1}{2l}+2\delta^{1/r}\right)  n^{2}\\
&  \leq\left(  \frac{\varepsilon^{5}}{8}+\frac{\varepsilon^{5}}{8}\right)
n^{2}\leq\frac{\varepsilon^{3}}{4}\left\vert V^{\prime\prime}\right\vert
^{2}<\varepsilon^{3}\binom{\left\vert V^{\prime\prime}\right\vert }{2}.
\end{align*}

On the other hand, (\ref{defs}) and (\ref{defl}) imply%
\begin{equation}
s\leq\frac{\varepsilon n}{4kL\left(  \varepsilon^{3},r\right)  }\leq
\frac{\varepsilon n}{4lL\left(  \varepsilon^{3},r\right)  }\leq\varepsilon
^{2}n<\varepsilon\left\vert V^{\prime\prime}\right\vert . \label{sup}%
\end{equation}
Hence, we apply the Scooping Lemma to the graph $G\left[  V^{\prime\prime
}\right]  ,$ and find a partition $V^{\prime\prime}=\cup_{i=0}^{g}X_{i}$ such
that $\left\vert X_{0}\right\vert <\left\lceil \varepsilon\left\vert
V^{\prime\prime}\right\vert \right\rceil ,$ and $\left\vert X_{1}\right\vert
=...=\left\vert X_{g}\right\vert =s.$ Set
\begin{align*}
q  &  =h+g,\\
W_{0}  &  =V_{0}\cup X_{0}\cup X^{\prime},\\
W_{h+i}  &  =X_{i},\text{ }i\in\left[  g\right]  .
\end{align*}
From (\ref{x1up}) and (\ref{sup}) it follows
\[
\left\vert W_{0}\right\vert =\left\vert V_{0}\right\vert +\left\vert
X_{0}\right\vert +\left\vert X^{\prime}\right\vert <M\left(  \delta,l\right)
+\left\lceil \varepsilon\left\vert V^{\prime\prime}\right\vert \right\rceil
+3\varepsilon n<6\varepsilon n.
\]
Finally, (\ref{props}) implies
\[
q=h+g\leq\frac{n}{s}<L\left(  \varepsilon,r+1\right)  ,
\]
completing the proof.
\end{proof}

\section{\label{Cor} Extensions of Theorem \ref{maint}}

Generally speaking, Theorem \ref{maint} states that, if certain conditions
about a graph are met, then its vertices can be partitioned in a specific way.
It turns out that, in addition, the partition may be selected to be
$\varepsilon$-uniform. This is is the topic of the following two theorems.

\begin{theorem}
\label{maint3} For all $\varepsilon>0$, $r\geq2$ and $k\geq2,$ there exist
$\rho=\rho\left(  \varepsilon,r,k\right)  >0$ and $K=K\left(  \varepsilon
,r,k\right)  $ such that, for every graph $G$ of sufficiently large order $n,$
the following assertion holds

If $k_{r}\left(  G\right)  \leq\rho n^{r},$ then there exists an $\varepsilon
$-uniform partition $V\left(  G\right)  =\cup_{i=0}^{q}V_{i}$ with $k\leq
q\leq K,$ and%
\[
e\left(  V_{i}\right)  \leq\varepsilon\binom{\left\vert V_{i}\right\vert }{2}%
\]
for every $i\in\left[  q\right]  .$
\end{theorem}

\begin{proof}
Our proof is essentially the same as the proof of Theorem \ref{maint}.

Suppose $M\left(  \varepsilon,l\right)  $ is as defined in SUL, and
$\xi\left(  \varepsilon,r\right)  ,$ $L\left(  \varepsilon,r\right)  $ are as
defined in Theorem \ref{maint}. Assume $\varepsilon$ sufficiently small and
set
\begin{align}
\delta &  =\min\left\{  \frac{\varepsilon^{2}}{8L\left(  \varepsilon
^{3},r\right)  },\frac{\varepsilon}{4}\right\} \label{defdel1}\\
l  &  =\max\left\{  k,\frac{2}{\varepsilon}\right\}  .\label{defl1}\\
K  &  =K\left(  \varepsilon,r,k\right)  =\frac{8M\left(  \delta,l\right)
L\left(  \varepsilon^{3},r\right)  }{\varepsilon}\label{defK}\\
\rho &  =\rho\left(  \varepsilon,r,k\right)  =\frac{\xi\left(  \varepsilon
^{3},r\right)  }{\left(  2M\left(  \delta,l\right)  \right)  ^{r}}.
\label{defro}%
\end{align}
Let $G$ be a graph of sufficiently large order $n,$ and let $k_{r}\left(
G\right)  \leq\rho n^{r}.$ It suffices to find a partition $V\left(  G\right)
=\cup_{i=0}^{q}W_{i}$ such that:

\emph{(i)} $k\leq q\leq K;$\ 

\emph{(ii)} $\left\vert W_{0}\right\vert <3\varepsilon n,$ $\left\vert
W_{1}\right\vert =...=\left\vert W_{q}\right\vert ;$

\emph{(iii)} for every $i\in\left[  q\right]  ,$ $e\left(  W_{i}\right)
<\varepsilon\left\vert W_{i}\right\vert ^{2};$

\emph{(iv)} at most $\varepsilon q^{2}$ pairs are not $\varepsilon$-uniform.

Applying SUL, we find a $\delta$-uniform partition $V\left(  G\right)
=\cup_{i=0}^{p}V_{i}$ with $l\leq p\leq M\left(  \delta,l\right)  $. Set
$t=\left\vert V_{1}\right\vert $ and observe that
\begin{equation}
\frac{n}{2p}\leq\left(  1-\delta\right)  \frac{n}{p}<t\leq\frac{n}{p}.
\label{tbnd}%
\end{equation}
For every $i\in\left[  p\right]  ,$ we have%
\[
k_{r}\left(  V_{i}\right)  \leq k_{r}\left(  G\right)  <\rho n^{r}=\frac
{\xi\left(  \varepsilon^{3},r\right)  }{\left(  2M\left(  \delta,l\right)
\right)  ^{r}}n^{r}\leq\xi\left(  \varepsilon^{3},r\right)  \left(  \frac
{n}{2p}\right)  ^{r}\leq\xi\left(  \varepsilon^{3},r\right)  t^{r}.
\]
Hence, for every $i\in\left[  p\right]  ,$ we apply Theorem \ref{maint}, and
find an equitable partition $V_{i}=\cup_{j=0}^{m_{i}}Y_{ij}$ with $m_{i}\leq
L\left(  \varepsilon^{3},r\right)  ,$ and
\begin{equation}
e\left(  Y_{ij}\right)  \leq\varepsilon^{3}\binom{\left\vert Y_{ij}\right\vert
}{2} \label{lowd}%
\end{equation}
for every $j\in\left[  m_{i}\right]  $. Set%
\begin{equation}
s=\left\lfloor \frac{\varepsilon n}{4pL\left(  \varepsilon^{3},r\right)
}\right\rfloor \label{defs1}%
\end{equation}
and observe that
\begin{equation}
\frac{n}{s}<\frac{8pL\left(  \varepsilon^{3},r\right)  }{\varepsilon}%
<\frac{8M\left(  \delta,l\right)  L\left(  \varepsilon^{3},r\right)
}{\varepsilon}=K\left(  \varepsilon,r,k\right)  . \label{props1}%
\end{equation}
Also, for every $i\in\left[  p\right]  $ and $j\in\left[  m_{i}\right]  ,$
(\ref{defs1}) and (\ref{tbnd}) imply%
\[
s\leq\frac{\varepsilon n}{4pL\left(  \varepsilon^{3},r\right)  }%
<\varepsilon\frac{t}{2L\left(  \varepsilon^{3},r\right)  }\leq\varepsilon
\frac{t}{2m_{i}}<\varepsilon\left\lfloor \frac{t}{m_{i}}\right\rfloor
=\varepsilon\left\vert Y_{ij}\right\vert .
\]
Hence, for every $i\in\left[  p\right]  $ and $j\in\left[  m_{i}\right]  ,$ in
view of (\ref{lowd}), we apply the Scooping Lemma to the graph $G\left[
Y_{ij}\right]  ,$ and find a partition $Y_{ij}=\cup_{h=0}^{p_{ij}}W_{ijh}$
with $\left\vert W_{ij0}\right\vert \leq\left\lceil \varepsilon\left\vert
Y_{ij}\right\vert \right\rceil $ such that $\left\vert W_{ijh}\right\vert =s,$
and\ $e\left(  W_{ijh}\right)  <\varepsilon\binom{s}{2}$ for every
$h\in\left[  p_{ij}\right]  .$ Setting
\[
W_{0}=V_{0}\cup\left(  \cup_{i\in I}Y_{i0}\right)  \cup\left(  \cup_{i\in
I}\cup_{j=1}^{m_{i}}W_{ij0}\right)  ,
\]
we obtain
\begin{align}
\left\vert W_{0}\right\vert  &  =\left\vert V_{0}\right\vert +\left\vert
\left(  \cup_{i\in I}Y_{i0}\right)  \right\vert +\left\vert \left(  \cup_{i\in
I}\cup_{j=1}^{m_{i}}W_{ij0}\right)  \right\vert <\left\vert V_{0}\right\vert
+\sum_{i=1}^{p}m_{i}+\sum_{i=1}^{p}\sum_{j=1}^{m_{i}}\left\lceil
\varepsilon\left\vert Y_{ij}\right\vert \right\rceil \nonumber\\
&  <M\left(  \delta,l\right)  +pL\left(  \varepsilon^{3},r\right)
+2\varepsilon\sum_{i=1}^{p}m_{i}\left\lfloor \frac{t}{m_{i}}\right\rfloor
<3\varepsilon n. \label{w0up}%
\end{align}

Denote by $q$ the number of the sets $W_{ijh}$ $\left(  i\in I^{\prime},\text{
}j\in\left[  m_{i}\right]  ,\text{ }h\in\left[  p_{ij}\right]  \right)  ,$ and
renumber them sequentially from $1$ to $q$. Clearly, from (\ref{w0up}) and
(\ref{props1}), we have
\begin{equation}
\frac{\left(  1-3\varepsilon\right)  n}{s}\leq q\leq\frac{n}{s}\leq K\left(
\varepsilon,r,k\right)  . \label{minq}%
\end{equation}

Let us check that the partition $V\left(  G\right)  =\cup_{i=0}^{q}W_{i}$
satisfies \emph{(i)-(iv)}$.$ For every $i\in\left[  p\right]  $, the cluster
$V_{i}$ contains at least one $W_{j}$ $\left(  j\in\left[  q\right]  \right)
$, so \emph{(i)} holds. Observe that $\left\vert W_{0}\right\vert
\,<3\varepsilon n,$ and
\[
\left\vert W_{i}\right\vert =s,\text{ \ \ }e\left(  W_{i}\right)
<\varepsilon\binom{\left\vert W_{i}\right\vert }{2}%
\]
for every $i\in\left[  q\right]  ,$ so \emph{(ii)} and \emph{(iii)} also hold.
To complete the proof, it remains to check \emph{(iv)}. Suppose $W_{a}\subset
V_{i},$ $W_{b}\subset V_{j}.$ If the pair $\left(  V_{i},V_{j}\right)  $ is
$\delta$-uniform, then,%
\[
\left\vert W_{a}\right\vert =\left\vert W_{b}\right\vert =s\geq\frac
{\varepsilon n}{8kL\left(  \varepsilon^{3},r\right)  }\geq\frac{\varepsilon
}{8L\left(  \varepsilon^{3},r\right)  }t.
\]
\ Since (\ref{defdel1}) implies
\[
\varepsilon=\max\left\{  \frac{8L\left(  \varepsilon^{3},r\right)
}{\varepsilon}\delta,2\delta\right\}  ,
\]
from Lemma \ref{Slicel}, it follows that the pair $\left(  W_{a},W_{b}\right)
$ is $\varepsilon$-uniform. Therefore, if the pair $\left(  W_{a}%
,W_{b}\right)  $ is not $\varepsilon$-uniform, then either $i=j$ or the pair
$\left(  V_{i},V_{j}\right)  $ is not $\delta$-uniform. For every $i\in\left[
p\right]  ,$ $V_{i}$ contains at most $\left\lfloor t/s\right\rfloor $ sets
$W_{a},$ so the number of the pairs $\left(  W_{a},W_{b}\right)  $ that are
not $\varepsilon$-uniform is at most%
\[
p\binom{\left\lfloor t/s\right\rfloor }{2}+\delta p^{2}\left\lfloor
t/s\right\rfloor ^{2}<\left(  \frac{1}{2p}+\delta\right)  \left(  \frac{pt}%
{s}\right)  ^{2}\leq\left(  \frac{1}{2l}+\delta\right)  \left(  \frac{n}%
{s}\right)  ^{2}.
\]
From (\ref{defdel1}), (\ref{defL1}) and (\ref{minq}) we find that
\[
\left(  \frac{1}{2l}+\delta\right)  \left(  \frac{n}{s}\right)  ^{2}\leq
\frac{\varepsilon}{2}\left(  \frac{n}{s}\right)  ^{2}<\varepsilon\left(
\frac{\left(  1-3\varepsilon\right)  n}{s}\right)  ^{2}\leq\varepsilon q^{2},
\]
completing the proof.
\end{proof}

Applying routine argument, we obtain the following corollary.

\begin{theorem}
\label{maint2} For all $\varepsilon>0$, $r\geq2$ and $k\geq2,$ there exist
$\rho=\rho\left(  \varepsilon,r,k\right)  >0$ and $K=K\left(  \varepsilon
,r,k\right)  $ such that, for every graph $G$ of sufficiently large order $n,$
the following assertion holds.

If $k_{r}\left(  G\right)  \leq\rho n^{r},$ then there exists a partition
$V\left(  G\right)  =\cup_{i=1}^{q}V_{i}$ with $k\leq q\leq K$ such that%
\[
\left\lfloor n/q\right\rfloor \leq\left\vert V_{i}\right\vert \leq\left\lceil
n/q\right\rceil ,\text{ \ \ \ \ }e\left(  V_{i}\right)  \leq\varepsilon
\binom{\left\vert V_{i}\right\vert }{2}%
\]
for every $i\in\left[  q\right]  ,$ and at most $\varepsilon q^{2}$ pairs
$\left(  W_{i},W_{j}\right)  $ are not $\varepsilon$-uniform.
\end{theorem}

\section{\label{Bipar} Bipartitions of low density}

In this section we shall deduces Theorem \ref{Ers1}. We first state and prove
a preliminary result of its own interest. In fact, this is a particular result
on judicious bipartitions of dense graphs with moderately many $r$-cliques; it
significantly differs from known general results as in \cite{BoSc99},
\cite{BoSc02} and \cite{ABKS03}.

\begin{theorem}
\label{Ers2} For all $r\geq3$, $c>0$ and $\varepsilon>0,$ there exist $\xi
=\xi(\varepsilon,c,r)>0$ and $\beta=\beta(\varepsilon,c,r)>0$ such that, for
every graph $G=G(n,\left\lfloor cn^{2}\right\rfloor )$ of sufficiently large
order $n,$ the following assertion holds.

If $k_{r}\left(  G\right)  <\xi n^{r}$, then there exists a partition
$V\left(  G\right)  =V_{1}\cup V_{2}$ such that
\[
e\left(  V_{1}\right)  <\varepsilon\left\vert V_{1}\right\vert ^{2}\text{
\ \ and \ \ }e\left(  V_{2}\right)  <\left(  c-\beta\right)  \left\vert
V_{2}\right\vert ^{2}.
\]

\end{theorem}

\begin{proof}
Let $\xi_{1}\left(  \varepsilon,r\right)  $ and $L_{1}\left(  \varepsilon
,r\right)  $ correspond to $\xi\left(  \varepsilon,r\right)  $ and $L\left(
\varepsilon,r\right)  $ as defined in Theorem \ref{maint}. Set%
\begin{align}
\sigma &  =\min\left\{  c/4,\varepsilon\right\} \nonumber\\
\xi &  =\xi\left(  \varepsilon,c,r\right)  =\xi_{1}\left(  \sigma,r\right)
\nonumber\\
L  &  =L\left(  \varepsilon,c,r\right)  =L_{1}\left(  \sigma,r\right)
\nonumber\\
\beta &  =\beta(\varepsilon,c,r)=\frac{c-2\sigma}{L^{2}} \label{defb}%
\end{align}

Let the graph $G=G\left(  n,\left\lfloor cn^{2}\right\rfloor \right)  $ be
with $k_{r}\left(  G\right)  <\xi n^{r}.\,\ $If $n$ is sufficiently large, we
apply Theorem \ref{maint3}, and find a $\sigma$-uniform partition $V\left(
G\right)  =\cup_{i=0}^{k}W_{i}$ with $k\leq L$ such that $\left\vert
W_{i}\right\vert =\left\lfloor n/k\right\rfloor $ and $e\left(  W_{i}\right)
\leq\sigma\left\lfloor n/k\right\rfloor ^{2}$ for every $i\in\left[  k\right]
$. Set $V^{\prime}=V\backslash W_{0}.$ We may and shall assume that the
cluster $W_{1}$ satisfies
\[
e\left(  W_{1},V^{\prime}\backslash W_{1}\right)  =\max_{i\in\left[  k\right]
}\left\{  e\left(  W_{i},V^{\prime}\backslash W_{i}\right)  \right\}  .
\]
We shall prove that the partition $V\left(  G\right)  =V_{1}\cup V_{2},$
defined with $V_{1}=W_{1},$ $V_{2}=V\backslash W_{1},$ satisfies the
requirements. Indeed, we immediately have%
\[
e\left(  V_{1}\right)  <\sigma\left\vert V_{1}\right\vert ^{2}\leq
\varepsilon\left\vert V_{1}\right\vert ^{2}.
\]
Therefore, all we have to prove is that, for $n$ sufficiently large, $e\left(
V_{2}\right)  <\left(  c-\beta\right)  \left\vert V_{2}\right\vert ^{2},$ that
is to say
\begin{equation}
e\left(  V\backslash W_{1}\right)  <\left(  c-\beta\right)  \left(
n-\left\lfloor \frac{n}{k}\right\rfloor \right)  ^{2}. \label{topr}%
\end{equation}

We have,
\begin{equation}
e\left(  V^{\prime}\right)  =e\left(  V\right)  -e\left(  W_{0},V^{\prime
}\right)  -e\left(  W_{0}\right)  \geq e\left(  V\right)  -\left\vert
W_{0}\right\vert n>e\left(  V\right)  -kn. \label{loev}%
\end{equation}

On the other hand,%
\begin{align*}
e\left(  V^{\prime}\right)   &  =\sum_{i=1}^{k}e\left(  W_{i}\right)
+\frac{1}{2}\sum_{i=1}^{k}e\left(  W_{i},V^{\prime}\backslash W_{i}\right)
\leq\sum_{i=1}^{k}e\left(  W_{i}\right)  +\frac{k}{2}e\left(  W_{1},V^{\prime
}\backslash W_{1}\right) \\
&  \leq\sigma k\left\lfloor \frac{n}{k}\right\rfloor ^{2}+\frac{k}{2}e\left(
W_{1},V^{\prime}\backslash W_{1}\right)  .
\end{align*}
Therefore,%
\[
e\left(  W_{1},V^{\prime}\backslash W_{1}\right)  \geq\frac{2e\left(
V^{\prime}\right)  }{k}-2\sigma\left\lfloor \frac{n}{k}\right\rfloor ^{2}.
\]
This, together with (\ref{loev}), implies
\begin{align*}
e\left(  V\backslash W_{1}\right)   &  =e\left(  V\right)  -e\left(
W_{1}\right)  -e\left(  W_{1},V\backslash W_{1}\right)  \leq e\left(
V\right)  -e\left(  W_{1},V^{\prime}\backslash W_{1}\right) \\
&  \leq e\left(  V\right)  -\frac{2e\left(  V^{\prime}\right)  }{k}%
+2\sigma\left\lfloor \frac{n}{k}\right\rfloor ^{2}\leq\frac{k-2}{k}e\left(
V\right)  +kn+2\sigma\left\lfloor \frac{n}{k}\right\rfloor ^{2}.
\end{align*}
Hence, in view of $e\left(  V\right)  =\left\lfloor cn^{2}\right\rfloor ,$ we
deduce%
\begin{equation}
e\left(  V\backslash W_{1}\right)  \leq\frac{k-2}{k}\left\lfloor
cn^{2}\right\rfloor +kn+2\sigma\left\lfloor \frac{n}{k}\right\rfloor ^{2}.
\label{upev}%
\end{equation}

Assume that there are arbitrary large values of $n$ for which (\ref{topr}) is
false, thus%
\[
e\left(  V\backslash W_{1}\right)  \geq\left(  c-\beta\right)  \left(
n-\left\lfloor \frac{n}{k}\right\rfloor \right)  ^{2}%
\]
holds. Hence, in view of (\ref{upev}), we find that%

\[
\left(  c-\beta\right)  \left(  n-\left\lfloor \frac{n}{k}\right\rfloor
\right)  ^{2}\leq\frac{k-2}{k}\left\lfloor cn^{2}\right\rfloor +kn+2\sigma
\left\lfloor \frac{n}{k}\right\rfloor ^{2}.
\]
Dividing both sides by $n^{2}$ and taking the limit, we deduce%
\[
\left(  c-\beta\right)  \left(  1-\frac{1}{k}\right)  ^{2}\leq\frac{k-2}%
{k}c+\frac{2\sigma}{k^{2}},
\]
and hence,
\[
c-2\sigma\leq\beta\left(  k-1\right)  ^{2}<\beta L^{2},
\]
a contradiction with (\ref{defb}). The proof is completed.
\end{proof}

In \cite{Nik01}, for every graph $G=G(n,m),$ the function
\[
\Phi\left(  G,k\right)  =\min_{U\subset V\left(  G\right)  ,\text{ }\left\vert
U\right\vert =k}\left\{  \frac{e\left(  U\right)  }{k}+\frac{e\left(
V\backslash U\right)  }{n-k}-\frac{m}{n}\right\}
\]
is introduced, and, it is shown that, if $1\leq k\leq\left\lfloor
n/2\right\rfloor ,$ then
\[
\Phi\left(  G,\left\lfloor \frac{n}{2}\right\rfloor \right)  \leq\frac{k}%
{n-k}\Phi\left(  G,k\right)  .
\]

This inequality, together with Theorem \ref{Ers2}, easily implies Theorem
\ref{Ers}.

Note that $k_{r}\left(  \overline{G}\right)  $ is exactly the number of
independent $r$-sets in $G.$ Restating the Scooping Lemma, and Theorems
\ref{maint}, \ref{Ers1}, \ref{maint3} and \ref{maint2}, for the complementary
graph, we obtain equivalent assertions for graphs with few independent
$r$-sets; for example, the following theorem is equivalent to Theorem
\ref{Ers1}.

\begin{theorem}
\label{Ers1_com} For all $c>0$ and $r\geq3,$ there exist $\xi=\xi\left(
c,r\right)  >0$ and $\beta=\beta(c,r)>0$ such that, for $n$ sufficiently
large, and every graph $G=G(n,m)$ with $m\geq cn^{2}$, the following assertion holds.

If $k_{r}\left(  \overline{G}\right)  <\xi n^{r}$, then there exists a
partition $V\left(  G\right)  =V_{1}\cup V_{2}$ with $\left\vert
V_{1}\right\vert =\left\lfloor n/2\right\rfloor ,$ $\left\vert V_{2}%
\right\vert =\left\lceil n/2\right\rceil ,$ and
\[
e\left(  V_{1},V_{2}\right)  <\left(  1/2-\beta\right)  m.
\]

\end{theorem}

\section{\label{Spar}Refining partitions}

This section contains a proof of Theorem \ref{Rams}.

\begin{proof}
[Proof of Theorem \ref{Rams}]We follow essentially the proof of Theorem
\ref{maint3}.

Assume $\xi\left(  \varepsilon,r\right)  $ and $L\left(  \varepsilon,r\right)
$ as defined in Theorem \ref{maint}; assume $\varepsilon$ sufficiently small
and set
\begin{align}
\delta &  =\delta\left(  \varepsilon,r\right)  =\min\left\{  \frac
{\varepsilon^{2}}{8L\left(  \varepsilon^{3},r\right)  },\frac{\varepsilon}%
{8}\right\} \label{defdel2}\\
\rho &  =\rho\left(  \varepsilon,r\right)  =\xi\left(  \varepsilon
^{3},r\right) \label{defro1}\\
K  &  =K\left(  \varepsilon,r,k\right)  =\frac{8kL\left(  \varepsilon
^{3},r\right)  }{\varepsilon} \label{defK1}%
\end{align}

Let $G$ be a graph of sufficiently large order $n,$ and let $V\left(
G\right)  =\cup_{i=0}^{k}V_{i}$ be a $\delta$-uniform partition such that%
\begin{equation}
k_{r}\left(  V_{i}\right)  \leq\rho\left\lfloor n/k\right\rfloor ^{r}\text{
\ \ or \ \ }k_{r}\left(  \overline{G\left[  V_{i}\right]  }\right)  \leq
\rho\left\lfloor n/k\right\rfloor ^{r} \label{monos}%
\end{equation}
for every $i\in\left[  k\right]  .$ We shall also assume that%
\begin{equation}
k>\frac{2}{\varepsilon}, \label{mink}%
\end{equation}
as, changing $\delta,$ $\rho$ and $K$ appropriately, we may refine the
partition $V\left(  G\right)  =\cup_{i=0}^{k}V_{i}$ so that (\ref{monos}) and
(\ref{mink}) hold. To prove the theorem, it suffices to find a partition
$V\left(  G\right)  =\cup_{i=0}^{q}W_{i}$ such that:

\emph{(i)} $k\leq q\leq K;$\ 

\emph{(ii)} $\left\vert W_{0}\right\vert <3\varepsilon n,$ $\left\vert
W_{1}\right\vert =...=\left\vert W_{q}\right\vert ;$

\emph{(iii)} for every $i\in\left[  q\right]  ,$ $e\left(  W_{i}\right)
<\varepsilon\binom{\left\vert W_{i}\right\vert }{2}$ or $e\left(
W_{i}\right)  <\left(  1-\varepsilon\right)  \binom{\left\vert W_{i}%
\right\vert }{2};$

\emph{(iv)} at most $\varepsilon q^{2}$ pairs are not $\varepsilon$-uniform.

Set $t=\left\vert V_{1}\right\vert $ and observe that
\begin{equation}
\frac{n}{2k}\leq\left(  1-\delta\right)  \frac{n}{k}<t\leq\frac{n}{k}.
\label{tbnd1}%
\end{equation}
From (\ref{defro1}) and (\ref{monos}) it follows
\[
k_{r}\left(  V_{i}\right)  \leq\xi\left(  \varepsilon^{3},r\right)
t^{r}\text{ \ \ or \ \ }k_{r}\left(  \overline{G\left[  V_{i}\right]
}\right)  \leq\xi\left(  \varepsilon^{3},r\right)  t^{r}%
\]
for every $i\in\left[  k\right]  $. Hence, for every $i\in\left[  k\right]  ,$
we apply Theorem \ref{maint} to the graph $G\left(  V_{i}\right)  $ or to its
complement, and find an equitable partition $V_{i}=\cup_{j=0}^{m_{i}}Y_{ij}$
with $m_{i}\leq L\left(  \varepsilon^{3},r\right)  $ such that
\[
e\left(  Y_{ij}\right)  \leq\varepsilon^{3}\binom{\left\vert Y_{ij}\right\vert
}{2}\text{ \ \ or \ \ }e\left(  Y_{ij}\right)  \geq\left(  1-\varepsilon
^{3}\right)  \binom{\left\vert Y_{ij}\right\vert }{2}%
\]
for every $j\in\left[  m_{i}\right]  $. Set%
\begin{equation}
s=\left\lfloor \frac{\varepsilon n}{4kL\left(  \varepsilon^{3},r\right)
}\right\rfloor \label{defs2}%
\end{equation}
and observe that
\begin{equation}
\frac{n}{s}<\frac{8kL\left(  \varepsilon^{3},r\right)  }{\varepsilon}=K\left(
\varepsilon,r,k\right)  . \label{props2}%
\end{equation}
Also, note that, for every $i\in\left[  k\right]  ,$ $j\in\left[
m_{i}\right]  ,$ (\ref{tbnd1}) and (\ref{defs2}) imply
\[
s\leq\frac{\varepsilon n}{4kL\left(  \varepsilon^{3},r\right)  }%
<\varepsilon\frac{t}{2K\left(  \varepsilon^{3},r,l\right)  }\leq
\varepsilon\frac{t}{2m_{i}}\leq\varepsilon\left\lfloor \frac{t}{m_{i}%
}\right\rfloor =\varepsilon\left\vert Y_{ij}\right\vert .
\]
Hence, for every $i\in\left[  k\right]  ,$ $j\in\left[  m_{i}\right]  ,$ we
apply the Scooping Lemma to the graph $G\left[  Y_{ij}\right]  $ or its
complement, and find a partition $Y_{ij}=\cup_{h=0}^{p_{ij}}W_{ijh}$ with
$\left\vert W_{ij0}\right\vert \leq\left\lceil \varepsilon\left\vert
Y_{ij}\right\vert \right\rceil $ such that
\[
\left\vert W_{ijh}\right\vert =s,\text{ \ \ and \ \ }e\left(  W_{ijh}\right)
<\varepsilon\binom{s}{2}\text{ \ \ or \ \ }e\left(  W_{ijh}\right)  >\left(
1-\varepsilon\right)  \binom{s}{2}%
\]
for every $h\in\left[  p_{ij}\right]  $. Setting
\[
W_{0}=V_{0}\cup\left(  \cup_{i=1}^{k}Y_{i0}\right)  \cup\left(  \cup_{i=1}%
^{k}\cup_{j=1}^{m_{i}}W_{ij0}\right)  ,
\]
we obtain
\begin{align}
\left\vert W_{0}\right\vert  &  =\left\vert V_{0}\right\vert +\left\vert
\left(  \cup_{i\in I}Y_{i0}\right)  \right\vert +\left\vert \left(  \cup_{i\in
I}\cup_{j=1}^{m_{i}}W_{ij0}\right)  \right\vert <\left\vert V_{0}\right\vert
+\sum_{i=1}^{k}m_{i}+\sum_{i=1}^{k}\sum_{j=1}^{m_{i}}\left\lceil
\varepsilon\left\vert Y_{ij}\right\vert \right\rceil \nonumber\\
&  <k+kL\left(  \varepsilon^{3},r\right)  +2\varepsilon\sum_{i=1}^{k}%
m_{i}\left\lfloor \frac{t}{m_{i}}\right\rfloor <3\varepsilon n. \label{w0up1}%
\end{align}

Denote by $q$ be the number of the sets $W_{ijh}$ $\left(  i\in I^{\prime
},\text{ }j\in\left[  m_{i}\right]  ,\text{ }h\in\left[  p_{ij}\right]
\right)  ,$ and renumber them sequentially from $1$ to $q$. Clearly, from
(\ref{w0up1}) and (\ref{props2}), we have
\begin{equation}
\frac{\left(  1-3\varepsilon\right)  n}{s}\leq q\leq\frac{n}{s}\leq K\left(
\varepsilon,r,k\right)  . \label{mind1}%
\end{equation}

Let us check that the partition $V\left(  G\right)  =\cup_{i=0}^{q}W_{i}$
satisfies \emph{(i)-(iv)}$.$ For every $i\in\left[  k\right]  $, the cluster
$V_{i}$ contains at least one $W_{j}$ $\left(  j\in\left[  q\right]  \right)
$, so \emph{(i)} holds. Observe that $\left\vert W_{0}\right\vert
<3\varepsilon n,$ and%
\[
\left\vert W_{i}\right\vert =s,\text{ \ \ \ }e\left(  W_{i}\right)
<\varepsilon\binom{s}{2}%
\]
for every $i\in\left[  q\right]  ,$ so \emph{(ii)} and \emph{(iii)} also hold.
To complete the proof, it remains to check \emph{(iv)}. Suppose $W_{a}\subset
V_{i},$ $W_{b}\subset V_{j}.$ If the pair $\left(  V_{i},V_{j}\right)  $ is
$\delta$-uniform then%
\[
\left\vert W_{a}\right\vert =\left\vert W_{b}\right\vert =s\geq\frac
{\varepsilon n}{8kL\left(  \varepsilon^{3},r\right)  }\geq\frac{\varepsilon
}{8L\left(  \varepsilon^{3},r\right)  }t.
\]
\ Since, from (\ref{defdel2}), we have
\[
\varepsilon=\max\left\{  \frac{8L\left(  \varepsilon^{3},r\right)
}{\varepsilon}\delta,2\delta\right\}  ,
\]
Lemma \ref{Slicel} implies that the pair $\left(  W_{a},W_{b}\right)  $ is
$\varepsilon$-uniform. Therefore, if the pair $\left(  W_{a},W_{b}\right)  $
is not $\varepsilon$-uniform, then either $i=j,$ or the pair $\left(
V_{i},V_{j}\right)  $ is not $\delta$-uniform. For every $i\in\left[
k\right]  $, $V_{i}$ contains at most $\left\lfloor t/s\right\rfloor $ sets
$W_{a},$ so the number of the pairs $\left(  W_{a},W_{b}\right)  $ that are
not $\varepsilon$-uniform is at most%
\[
k\binom{\left\lfloor t/s\right\rfloor }{2}+\delta k^{2}\left\lfloor
t/s\right\rfloor ^{2}<\left(  \frac{1}{2k}+\delta\right)  \left(  \frac{kt}%
{s}\right)  ^{2}\leq\left(  \frac{1}{2k}+\delta\right)  \left(  \frac{n}%
{s}\right)  ^{2}.
\]
From (\ref{defdel2}) and (\ref{mink}) we find that
\[
\left(  \frac{1}{2k}+\delta\right)  \left(  \frac{n}{s}\right)  ^{2}\leq
\frac{\varepsilon}{2}\left(  \frac{n}{s}\right)  ^{2}<\varepsilon\left(
\frac{\left(  1-3\varepsilon\right)  n}{s}\right)  ^{2}\leq\varepsilon q^{2},
\]
completing the proof.
\end{proof}

\section{\label{Mxproof}\ Induced subgraphs}

In this section we shall prove Theorem \ref{maintx}. We start by a simple
partitioning lemma.

\begin{lemma}
\label{partL1}For all $\varepsilon>0$ and $b\geq2,$ there exist $\gamma
=\gamma\left(  \varepsilon,b\right)  $ and $n\left(  \varepsilon,b\right)  $
such that, for $n>n\left(  \varepsilon,b\right)  ,$ if the edges of $K_{n}$
are colored in red, blue and green, then the following assertion holds.

If there are fewer then $\gamma n^{2}$ green edges, then there exists a
partition $V\left(  K_{n}\right)  =\cup_{i=0}^{q}V_{i}$ such that $\left\vert
V_{0}\right\vert <\varepsilon n,$ $\left\vert V_{1}\right\vert =...=\left\vert
V_{q}\right\vert =b,$ and $V_{i}$ spans either a red or a blue $b$-clique for
every $i\in\left[  n\right]  $.
\end{lemma}

\begin{proof}
Ramsey's theorem implies that, for every $b,$ there exists $r=r\left(
b\right)  $ such that, if $n>r$ and the edges of $K_{n}$ are colored in two
colors, then there exists a monochromatic $K_{b}$.

We shall assume $\varepsilon<1,$ else there is nothing to prove. Set%
\[
\gamma=\gamma\left(  \varepsilon,b\right)  =\frac{\varepsilon^{2}}{4r}.
\]
Suppose $n>2r/\varepsilon^{2}$ and let the edges of $K_{n}$ be colored in red,
blue and green, so that there are fewer then $\gamma n^{2}$ green edges.
Therefore, there are at least
\[
\binom{n}{2}-\gamma n^{2}=\frac{n^{2}}{2}-\frac{n}{2}-\frac{n^{2}}{4r}%
>\frac{\left(  r-1\right)  n^{2}}{2r}%
\]
red or blue edges. Hence, by Tur\'{a}n's theorem, there is a set $U$ of
cardinality $r+1$ inducing only red or blue edges. By the choice of $r,$ $U$
induces a red or a blue $b$-clique; select one and denote its vertex set by
$V_{1}$. Proceed selecting sets $V_{2},...,V_{q}$ as follows: having selected
$V_{1},...,V_{i},$ if%
\[
\binom{n-bi}{2}-\frac{\varepsilon^{2}}{4r}n^{2}<\frac{\left(  r-1\right)
\left(  n-bi\right)  ^{2}}{2r}%
\]
stop the sequence, else, by Tur\'{a}n's theorem, find a set $U$ of cardinality
$r+1$ inducing only red or blue edges. By the choice of $r,$ $U$ induces a red
or a blue $b$-clique; select one and denote its vertex set by $V_{i+1}$.

Let $V_{q}$ be the last selected set; set $V_{0}=V\left(  G\right)
\backslash\left(  \cup_{i=1}^{q}V_{i}\right)  .$ The stop condition implies
\[
\frac{\left(  r-1\right)  \left\vert V_{0}\right\vert ^{2}}{2r}\geq
\binom{\left\vert V_{0}\right\vert }{2}-\frac{\varepsilon^{2}}{4r}n^{2}%
>\frac{\left\vert V_{0}\right\vert ^{2}}{2}-\frac{n}{2}-\frac{\varepsilon^{2}%
}{4r}n^{2}.
\]
This, and $\varepsilon^{2}n<2r,$ imply $\left\vert V_{0}\right\vert
<\varepsilon n$. Every set $V_{1},...,V_{q}$ spans either a red or a blue
$b$-clique, so the partition $V\left(  G\right)  =\cup_{i=0}^{k}V_{i}$ is as required.
\end{proof}

We shall need also the following modification of Lemma \ref{XPle}.

\begin{lemma}
\label{XPle1} Let $r\geq1,$ $0<2\varepsilon^{1/r}<d<1-2\varepsilon^{1/r},$ and
let $\left(  A,B\right)  $ be an $\varepsilon$-uniform pair with $d\left(
A,B\right)  =d.$ There are at most $\varepsilon2^{r}\left\vert A\right\vert
^{r}$ $r$-sets $R\subset A$ such that, there exists a partition $R=R_{0}\cup
R_{1}$ satisfying%
\[
\left\vert \left(  \cap_{u\in R_{0}}\Gamma\left(  u\right)  \right)
\cap\left(  \cap_{u\in R_{1}}\left(  B\backslash\Gamma\left(  u\right)
\right)  \right)  \cap B\right\vert \leq\varepsilon\left\vert B\right\vert .
\]

\end{lemma}

\begin{proof}
[Proof of Theorem \ref{maintx}]For $r=2$ the assertion easily follows from the
Scooping Lemma. To prove it for $r>2$ we apply induction on $r$ - assuming it
holds for $r$, we shall prove it for $r+1.$

It is sufficient to find $\xi=\xi\left(  \varepsilon,r+1\right)  >0$ and
$L=L\left(  \varepsilon,r+1\right)  $ such that, for every graph $H$ of order
$r+1,$ and every graph $G$ of sufficiently large order $n,$ if $k_{H}\left(
G\right)  <\xi n^{r+1},$ then there exists a partition $V\left(  G\right)
=\cup_{i=0}^{q}W_{i}$ such that:

\emph{(i)} $q\leq L;$\ 

\emph{(ii)} $\left\vert W_{0}\right\vert <6\varepsilon n,$ $\left\vert
W_{1}\right\vert =...=\left\vert W_{q}\right\vert ;$

\emph{(iii)} for every $i\in\left[  q\right]  ,$ $e\left(  W_{i}\right)
<\varepsilon\binom{\left\vert W_{i}\right\vert }{2}$ or $e\left(
W_{i}\right)  >\left(  1-\varepsilon\right)  \binom{\left\vert W_{i}%
\right\vert }{2}.$

We shall outline our proof first. Choose $b,$ $\delta,$ $\xi\left(
\varepsilon,r+1\right)  $, and $L\left(  \varepsilon,r+1\right)  $
appropriately. Select any graph $H,$ let $G$ be a graph of sufficiently large
order $n,$ and let $k_{H}\left(  G\right)  <\xi\left(  \varepsilon,r+1\right)
n^{r+1}.$ Applying SUL, we find a $\delta$-uniform partition $V\left(
G\right)  =\cup_{i=0}^{k}V_{i}.$ Note that, if $k_{F}\left(  V_{i}\right)  $
is proportional to $n^{r}$, and $V_{i}$ is incident to a $\delta$-uniform pair
of medium density, then, by Lemma \ref{XPle1}, there are substantially many
induced copies of $H$ in $G.$ Therefore, for every $V_{i},$ either
$k_{F}\left(  V_{i}\right)  $ is small, or $V_{i}$ is incident only to very
sparse or very dense $\delta$-uniform pairs.

Let $V^{\prime}$ be the vertices in the clusters $V_{i}$ with small
$k_{F}\left(  V_{i}\right)  ;$ set $V^{\prime\prime}=V\backslash\left(
V^{\prime}\cup V_{0}\right)  $.

\emph{1 Partitioning of }$V^{\prime}$

By the induction hypothesis, we partition each $V_{i}$ with small
$k_{F}\left(  V_{i}\right)  $ into a bounded number of sets $Y_{ij}$ that are
either very sparse or very dense and a small exceptional set; the exceptional
sets are collected in $X^{\prime}.$ Although the sets $Y_{ij}$ are very sparse
or very dense, they are not good for our purposes, for their cardinality may
vary with $i$. To overcome this obstacle, we first select a sufficiently small
integer $s$ proportional to $n.$ Then, by the Scooping Lemma, we partition
each of the sets $Y_{ij}$ into sparse or dense sets of cardinality exactly $s$
and a small exceptional set; the exceptional sets are added to $X^{\prime}$.

\emph{2 Partitioning of }$V^{\prime\prime}$

We partition $V^{\prime\prime}$ into dense or sparse sets of size $s$ and a
small exceptional set $X_{0}$. Suppose $V_{1},...,V_{g}$ are the clusters
whose union is $V^{\prime\prime}.$ If $g\ll k$, we let $X_{0}=V^{\prime\prime
},$ and complete the partition, so suppose that $g$ is proportional to $k.$
The density of the clusters is unknown, so we assemble them into larger groups
of $b$ clusters. Recall that the pairs $\left(  V_{i},V_{j}\right)  $ $\left(
1\leq i<j\leq g\right)  $ are either very sparse, or very dense, or are not
$\delta$-uniform. Color correspondingly the edges of $K_{g}$ in red, blue and
green. Since, there are fewer then $\delta k^{2}\leq\delta^{\prime}g^{2}$
green edges, by Lemma \ref{partL1}, we assemble almost all clusters
$V_{1},...,V_{g}$ into groups of exactly $b$ clusters and collect the vertices
of the few remaining clusters in a set $X_{0}.$ Observe that the pairs within
the same group are all either very dense or all very sparse. Finally, we apply
the Scooping Lemma to partition each of the groups into dense or sparse sets
of cardinality $s$ and an exceptional class; the exceptional classes are added
to $X_{0}$.

Let $W_{1},...,W_{q}$ be the sets of cardinality $s$ obtained during the
partitioning of $V^{\prime}$ and $V^{\prime\prime}$. Setting $W_{0}=V_{0}\cup
X_{0}\cup X^{\prime},$ the choice of $\delta$ implies $\left\vert
W_{0}\right\vert <\varepsilon n,$ so the partition $V\left(  G\right)
=\cup_{i=0}^{q}W_{i}$ satisfies \emph{(i)-(iii)}.

Let us give the details now. Let $\gamma\left(  \varepsilon,b\right)  $ and
$n\left(  \varepsilon,b\right)  $ be as defined in Lemma \ref{partL1},
$M\left(  \delta,l\right)  $ as defined in SUL, and $\xi\left(  \varepsilon
^{3},r\right)  ,$ $L\left(  \varepsilon^{3},r\right)  $ as defined in Theorem
\ref{maint}. Assume $\varepsilon$ sufficiently small and set%
\begin{align}
b  &  =\left\lceil \varepsilon^{-3}\right\rceil ,\label{defb1}\\
l  &  =n\left(  \varepsilon,b\right)  ,\label{defl3}\\
\delta &  =\min\left\{  \gamma\left(  \varepsilon,b\right)  \varepsilon
^{2},\frac{\xi\left(  \varepsilon^{3},r\right)  }{2^{r}+1},\frac
{\varepsilon^{3r}}{4^{r}}\right\} \label{del}\\
L\left(  \varepsilon,r+1\right)   &  =\frac{8M\left(  \delta,l\right)
L\left(  \varepsilon^{3},r\right)  }{\varepsilon},\label{defL2}\\
\xi &  =\xi\left(  \varepsilon,r+1\right)  =\frac{\delta^{2}}{\left(
2M\left(  \delta,l\right)  \right)  ^{r+1}}. \label{defksi1}%
\end{align}

Select any graph $H$ with $\left\vert H\right\vert =r+1$ and fix a vertex
$v\in V\left(  H\right)  .$ Set $F=H-v$ and let $F_{1}=\Gamma\left(  v\right)
,$ $F_{0}=F\backslash F_{1}.$

Let $G$ be a graph of sufficiently large order $n,$ and let $k_{H}\left(
G\right)  <\xi n^{r+1}.$ We apply SUL, and find a $\delta$-uniform partition
$V\left(  G\right)  =\cup_{i=0}^{k}V_{i}$ with $l\leq k\leq M\left(
\delta,l\right)  $. Set $t=\left\vert V_{1}\right\vert $ and observe that
\begin{equation}
\frac{n}{2k}\leq\left(  1-\delta\right)  \frac{n}{k}<t\leq\frac{n}{k}.
\label{bndt1}%
\end{equation}
Assume that there exist a cluster $V_{i}$ with $k_{F}\left(  V_{i}\right)
>\xi\left(  \varepsilon^{3},r\right)  t^{r},$ and a $\delta$-uniform pair
$\left(  V_{i},V_{j}\right)  $ with
\[
1-2\delta^{1/r}>d\left(  V_{i},V_{j}\right)  >2\delta^{1/r}.
\]
We apply Lemma \ref{XPle1} with $A=V_{i}$ and $B=V_{j},$ and find that there
are at least
\[
\xi\left(  \varepsilon^{3},r\right)  t^{r}-\delta2^{r}t^{r}\geq\delta t^{r}%
\]
induced subgraphs $V_{i}$ isomorphic to $F$ such that, if $X\subset G\left[
V_{i}\right]  $ and $\Phi:F\rightarrow X$ is an isomorphism, then
\[
\left\vert \left(  \cap_{u\in\Phi\left(  F_{0}\right)  }\Gamma\left(
u\right)  \right)  \cap\left(  \cap_{u\in\Phi\left(  F_{1}\right)  }\left(
B\backslash\Gamma\left(  u\right)  \right)  \right)  \cap B\right\vert
>\delta\left\vert B\right\vert .
\]
Hence, there are at least $\delta^{2}t^{r+1}$ induced copies of $H$ inducing a
copy of $F$ in $V_{i}$ and having a vertex in $V_{j}$. Therefore, from
(\ref{bndt1}) and (\ref{defksi1}), we find that
\[
k_{H}\left(  G\right)  \geq\delta^{2}t^{r+1}>\delta^{2}\left(  \frac{1-\delta
}{k}\right)  ^{r+1}n^{r+1}>\frac{\delta^{2}n^{r+1}}{\left(  2M\left(
\delta,l\right)  \right)  ^{r+1}}=\xi\left(  \varepsilon,r+1\right)  n^{r+1},
\]
a contradiction. Therefore, if $k_{F}\left(  V_{i}\right)  >\xi\left(
\varepsilon^{3},r\right)  t^{r},$ then every $\delta$-uniform pair $\left(
V_{i},V_{j}\right)  $ satisfies
\[
d\left(  V_{i},V_{j}\right)  \leq2\delta^{1/r}\text{ \ \ or \ \ }d\left(
V_{i},V_{j}\right)  \geq1-2\delta^{1/r}%
\]

Let
\[
I^{\prime}=\left\{  i:i\in\left[  k\right]  \text{, \ }k_{F}\left(
V_{i}\right)  \leq\xi\left(  \varepsilon^{3},r\right)  t^{r}\right\}  ,\text{
\ \ }I^{\prime\prime}=\left[  k\right]  \backslash I^{\prime}.
\]

First we shall partition the set $V^{\prime}=\cup_{i\in I^{\prime}}V_{i}.$ Set%
\begin{equation}
s=\left\lfloor \frac{\varepsilon n}{4kL\left(  \varepsilon^{3},r\right)
}\right\rfloor , \label{defs3}%
\end{equation}
and observe that
\begin{equation}
\frac{n}{s}<\frac{8kL\left(  \varepsilon^{3},r\right)  }{\varepsilon}%
<\frac{8M\left(  \delta,l\right)  L\left(  \varepsilon^{3},r\right)
}{\varepsilon}=L\left(  \varepsilon,r+1\right)  . \label{props3}%
\end{equation}
For every $i\in I^{\prime},$ by the induction hypothesis, we find an equitable
partition $V_{i}=\cup_{j=0}^{m_{i}}Y_{ij}$ with $\left\vert Y_{i0}\right\vert
<m_{i}\leq L\left(  \varepsilon^{3},r\right)  $ such that%
\[
e\left(  Y_{ij}\right)  <\varepsilon^{3}\binom{\left\vert Y_{ij}\right\vert
}{2}\text{ \ \ or \ \ }e\left(  Y_{ij}\right)  >\left(  1-\varepsilon
^{3}\right)  \binom{\left\vert Y_{ij}\right\vert }{2}%
\]
for every $j\in\left[  m_{i}\right]  .$ Also, for every $i\in I^{\prime}$ and
$j\in\left[  m_{i}\right]  ,$ (\ref{defs3}) and (\ref{bndt1}) imply%
\[
s\leq\frac{\varepsilon n}{4kL\left(  \varepsilon^{3},r\right)  }%
<\varepsilon\frac{t}{2L\left(  \varepsilon^{3},r\right)  }\leq\varepsilon
\frac{t}{2m_{i}}\leq\varepsilon\left\lfloor \frac{t}{m_{i}}\right\rfloor
=\varepsilon\left\vert Y_{ij}\right\vert .
\]
Hence, for every $i\in I^{\prime}$ and $j\in\left[  m_{i}\right]  ,$ we apply
the Scooping Lemma to the graph $G\left[  Y_{ij}\right]  ,$ and find a
partition $Y_{ij}=\cup_{q=0}^{p_{ij}}W_{ijq}$ with $\left\vert W_{ij0}%
\right\vert <\left\lceil \varepsilon\left\vert Y_{ij}\right\vert \right\rceil
$ such that
\[
\left\vert W_{ijq}\right\vert =s,\text{ \ \ and \ \ }e\left(  W_{ijq}\right)
<\varepsilon\binom{s}{2}\text{ \ \ or \ \ }e\left(  W_{ijq}\right)  >\left(
1-\varepsilon\right)  \binom{s}{2}%
\]
for every $q\in\left[  p_{ij}\right]  .$ Setting%
\[
X^{\prime}=\left(  \cup_{i\in I}Y_{i0}\right)  \cup\left(  \cup_{i\in I}%
\cup_{j=1}^{m_{i}}W_{ij0}\right)  ,
\]
we obtain
\begin{align}
\left\vert X^{\prime}\right\vert  &  =\left\vert \left(  \cup_{i\in I}%
Y_{i0}\right)  \right\vert +\left\vert \left(  \cup_{i\in I}\cup_{j=1}^{m_{i}%
}W_{ij0}\right)  \right\vert <\sum_{i\in I^{\prime}}m_{i}+\sum_{i\in
I^{\prime}}\sum_{j=1}^{m_{i}}\left\lceil \varepsilon\left\vert Y_{ij}%
\right\vert \right\rceil \nonumber\\
&  <kL\left(  \varepsilon^{3},r\right)  +2\varepsilon\sum_{i\in I^{\prime}%
}m_{i}\left\lfloor \frac{t}{m_{i}}\right\rfloor \leq kL\left(  \varepsilon
^{3},r\right)  +2\varepsilon n<3\varepsilon n. \label{x1up1}%
\end{align}

Denote by $h$ be the number of the sets $W_{ijq}$ $\left(  i\in I^{\prime
},\text{ }j\in\left[  m_{i}\right]  ,\text{ }q\in\left[  p_{ij}\right]
\right)  ,$ and renumber them sequentially from $1$ to $h$. \ Thus we have a
partition $V^{\prime}=X^{\prime}\cup\left(  \cup_{i=1}^{h}W_{i}\right)  $ with
$\left\vert X^{\prime}\right\vert \leq2\varepsilon n$ such that
\[
\left\vert W_{i}\right\vert =s,\ \ \text{and }\ \ e\left(  W_{i}\right)
<\varepsilon\binom{s}{2}\ \ \text{\ or }\ \ e\left(  W_{i}\right)  >\left(
1-\varepsilon\right)  \binom{s}{2}%
\]
for every $i\in\left[  h\right]  .$

Next we shall partition the set $V^{\prime\prime}=\cup_{i\in I^{\prime\prime}%
}V_{i}.$ For convenience assume $I^{\prime\prime}=\left[  g\right]  ;$ we may
assume $g\geq\varepsilon k,$ else, setting $W_{0}=V_{0}\cup X^{\prime}\cup
V^{\prime\prime},$ from (\ref{x1up1}), we have $W_{0}<4\varepsilon n,$\ and,
in view of (\ref{props3}), the proof is completed.

Recall that the pairs $\left(  V_{i},V_{j}\right)  $ $\left(  1\leq i<j\leq
g\right)  $ satisfy one of the following conditions:

a) $\left(  V_{i},V_{j}\right)  $ is $\delta$-uniform and $d\left(
V_{i},V_{i}\right)  <2\delta^{1/r};$

b) $\left(  V_{i},V_{j}\right)  $ is $\delta$-uniform and $d\left(
V_{i},V_{i}\right)  >1-2\delta^{1/r};$

c) $\left(  V_{i},V_{j}\right)  $ is not $\delta$-uniform.

Let $K_{g}$ be the complete graph on the vertex set $\left[  g\right]  ;$ for
every $1\leq i<j\leq g,$ color the edge $\left(  i,j\right)  $ in red, blue or
green correspondingly to a), b) and c). Observe that all pairs $\left(
V_{i},V_{j}\right)  $ $\left(  1\leq i<j\leq k\right)  $ that are not $\delta
$-uniform are fewer than%
\[
\delta k^{2}<\frac{\delta}{\varepsilon^{2}}g^{2}\leq\frac{\delta}%
{\varepsilon^{2}}g^{2}\leq\gamma\left(  \varepsilon,b\right)  g^{2},
\]
so, the green edges are fewer than $\gamma\left(  \varepsilon,b\right)
g^{2}.$ We apply Lemma \ref{partL1}, and find a partition $\left[  g\right]
=\cup_{i=0}^{a}X_{i}$ with $\left\vert X_{0}\right\vert <\varepsilon g,$ and
$X_{i}$ is either a red or a blue $b$-clique for every $i\in\left[  a\right]
$. For every $j=0,1,...,a,$ set $Y_{j}=\cup_{i\in X_{j}}V_{i}$; thus%
\begin{equation}
\left\vert Y_{0}\right\vert <\varepsilon gt\leq\varepsilon n,\text{
\ \ }\left\vert Y_{1}\right\vert =...=\left\vert Y_{a}\right\vert =bt.
\label{sizy}%
\end{equation}
Fix some $c\in\left[  a\right]  $ and assume $X_{c}$ a red $b$-clique. This is
to say that all pairs $\left(  V_{i},V_{j}\right)  ,$ $\left(  i,j\in
X_{c},\text{ }i<j\right)  ,$ are $\delta$-uniform and $d\left(  V_{i}%
,V_{i}\right)  <2\delta^{1/r}.$ Hence, from $\left\vert Y_{c}\right\vert =bt,$
(\ref{defb1}), and (\ref{del}), we deduce%
\begin{align*}
e\left(  Y_{c}\right)   &  =\sum_{i\in X_{c}}e\left(  V_{i}\right)
+\sum_{i,j\in X_{c},\text{ }i<j}e\left(  V_{i},V_{j}\right)  \leq b\binom
{t}{2}+\binom{b}{2}2\delta^{1/r}t^{2}\\
&  <\frac{bt^{2}}{2}+\delta^{1/r}b^{2}t^{2}=\left(  \frac{1}{2b}+\delta
^{1/r}\right)  \left\vert Y_{c}\right\vert ^{2}\leq\left(  \frac
{\varepsilon^{3}}{4}+\frac{\varepsilon^{3}}{4}\right)  \left\vert
Y_{c}\right\vert ^{2}=\varepsilon^{3}\binom{\left\vert Y_{c}\right\vert }{2}.
\end{align*}

On the other hand, (\ref{defs3}) and (\ref{sizy}) imply
\[
s\leq\frac{\varepsilon n}{4kL\left(  \varepsilon^{3},r\right)  }%
=\frac{\varepsilon}{4L\left(  \varepsilon^{3},r\right)  }t<\varepsilon
bt=\varepsilon\left\vert Y_{c}\right\vert .
\]
Hence, we apply the Scooping Lemma to the graph $G\left[  Y_{c}\right]  ,$ and
find a partition $Y_{c}=\cup_{i=0}^{f_{c}}Z_{ci}$ with $\left\vert
Z_{c0}\right\vert <s$ such that
\[
\left\vert Z_{ci}\right\vert =s,\text{ \ \ and \ \ }e\left(  Z_{ci}\right)
<\varepsilon\binom{s}{2}%
\]
for every $i\in\left[  f_{c}\right]  $.

If $X_{c}$ is a blue $b$-clique, proceeding in a similar way, we find a
partition $Y_{c}=\cup_{i=0}^{f_{c}}Z_{ci}$ with $\left\vert Z_{c0}\right\vert
<s$ such that
\[
\left\vert Z_{ci}\right\vert =s,\text{ \ \ and \ \ }e\left(  Z_{ci}\right)
>\left(  1-\varepsilon\right)  \binom{s}{2}%
\]
for every $i\in\left[  f_{c}\right]  $. Set $X^{\prime}=Y_{0}\cup\left(
\cup_{i=1}^{a}Z_{i0}\right)  $ and observe that (\ref{defs3}) and (\ref{sizy})
imply
\begin{align}
\left\vert X^{\prime}\right\vert  &  \leq\left\vert Y_{0}\right\vert
+\sum_{i=1}^{a}\left\vert Z_{i0}\right\vert \leq\varepsilon n+sa<\varepsilon
n+\frac{\varepsilon n}{4kL\left(  \varepsilon^{3},r\right)  }\frac{k}%
{b}\nonumber\\
&  <\left(  \varepsilon+\frac{\varepsilon^{4}}{4L\left(  \varepsilon
^{3},r\right)  }\right)  n<2\varepsilon n. \label{sup1}%
\end{align}

Denote by $h^{\prime}$ the number of the sets $Z_{ci}$ $\left(  c\in\left[
a\right]  ,\text{ }i\in\left[  f_{c}\right]  \right)  ,$ and renumber them
from $1$ to $h^{\prime}.$ Set
\begin{align*}
q  &  =h+h^{\prime},\\
W_{0}  &  =V_{0}\cup X_{0}\cup X^{\prime},\\
W_{h+i}  &  =Z_{i},\text{ }i\in\left[  h^{\prime}\right]  .
\end{align*}
From (\ref{x1up1}) and (\ref{sup1}) it follows
\[
\left\vert W_{0}\right\vert =\left\vert V_{0}\right\vert +\left\vert
X_{0}\right\vert +\left\vert X^{\prime}\right\vert <M\left(  \delta,l\right)
+2\varepsilon n+3\varepsilon n<6\varepsilon n.
\]
Finally, (\ref{props3}) implies
\[
q=h+h^{\prime}\leq\frac{n}{s}<L\left(  \varepsilon,r+1\right)  ,
\]
completing the proof
\end{proof}

Applying the same argument as in the proof of Theorem \ref{maint3}, we enhance
Theorem \ref{maintx} as follows.

\begin{theorem}
For all $\varepsilon>0,$ $r\geq2$ and $k\geq2,$ there exist $\rho=\rho\left(
\varepsilon,r,k\right)  >0$ and $K=K\left(  \varepsilon,r,k\right)  $ such
that, for every graph $G$ of sufficiently large order $n,$ the following
assertion holds$.$

If $k_{r}\left(  G\right)  \leq\rho n^{r},$ then there exists an $\varepsilon
$-uniform partition $V\left(  G\right)  =\cup_{i=0}^{q}V_{i}$ with $k\leq
q\leq K,$ and%
\[
e\left(  V_{i}\right)  <\varepsilon\binom{\left\vert V_{i}\right\vert }%
{2}\ \ \ \text{or}\ \ \ e\left(  V_{i}\right)  >\left(  1-\varepsilon\right)
\binom{\left\vert V_{i}\right\vert }{2}%
\]
for every $i\in\left[  q\right]  .$
\end{theorem}

\textbf{Acknowledgement.} The author is grateful to Cecil Rousseau for his
kind attention, and especially to B\'{e}la Bollob\'{a}s for his valuable
suggestions and remarks.


\begin{thebibliography}{99}                                                                                               %


\bibitem {Alo96}N. Alon, Bipartite subgraphs, \emph{Combinatorica} \textbf{16}
(1996), 301--311.

\bibitem {ABKS03}N. Alon, B. Bollob\'{a}s, M. Krivelevich and B. Sudakov,
Maximum cuts and judicious partitions in graphs without short cycles, \emph{J.
Comb. Theory Ser. B} \textbf{88} (2003), 329--346.

\bibitem {AKS03}N. Alon, M. Krivelevich and B. Sudakov, MaxCut in H-free
graphs, preprint.

\bibitem {Bol78}B. Bollob\'{a}s, \emph{Extremal Graph Theory,} Academic Press
Inc., London-New York, 1978, xx+488 pp.

\bibitem {Bol98}B. Bollob\'{a}s, \emph{Modern graph theory,} Graduate Texts in
Mathematics, \textbf{184,} Springer-Verlag, New York (1998), xiv+394 pp.

\bibitem {BoSc99}B. Bollob\'{a}s and A. D. Scott, Exact bounds for judicious
partitions of graphs, \emph{Combinatorica} \textbf{19} (1999), 473--486.

\bibitem {BoSc02}B. Bollob\'{a}s and A. D. Scott, Problems and results on
judicious partitions, \emph{Rand. Struct. Alg.} \textbf{21} (2002), 414--430.

\bibitem {Erd79}P. Erd\H{o}s, Some old and new problems in various branches of
combinatorics, \emph{Proceedings of the Tenth Southeastern Conference on
Combinatorics, Graph Theory and Computing }(Florida Atlantic Univ., Boca
Raton, Fla., 1979), \emph{Congress. Numer.} \textbf{XXIII--XXIV}, Utilitas
Math., Winnipeg, Man., 1979, pp. 19--37.

\bibitem {KoSi93}J. Koml\'{o}s and M. Simonovits, Szemer\'{e}di's regularity
lemma and its applications in graph theory, \emph{Combinatorics, Paul
Erd\H{o}s is Eighty}, Vol. \textbf{2} (Keszthely, 1993), pp. 295--352, Bolyai
Soc. Math. Stud., 2, \ J\'{a}nos Bolyai Math. Soc., Budapest, 1996.

\bibitem {KoRo03}Y. Kohayakawa and V. R\"{o}dl, Szemer\'{e}di's regularity
lemma and quasi-randomness, \emph{Recent advances in algorithms and
combinatorics}, CMS Books Math./Ouvrages Math. SMC, \textbf{11}, Springer, New
York, 2003, pp. 289--351.

\bibitem {Nik01}V. Nikiforov, On the edge distribution of a graph, \emph{Comb.
Prob. Comp. }\textbf{10} (2001), 543-555.
\end{thebibliography}
\end{document}